\documentclass[10pt]{article}
\usepackage{amsmath}
\usepackage[all,matrix,arrow,curve,poly,arc]{xy}

\newcommand{\qed}{%
  \ifmmode 
   \eqno{\qedsymbol}
  \else
    \leavevmode\unskip\penalty9999 \hbox{}\nobreak\hfill\hbox{\qedsymbol}
  \fi
}
\newcommand{\qedsymbol}{\leavevmode\vrule height 1.2ex width 1.1ex depth -.1ex}
\newenvironment{proof}{\begin{trivlist}\item[\hskip%
\labelsep{\bf Proof.\quad}]}%
{\hfill\qed\rm\end{trivlist}}

\newtheorem{theorem}{Theorem}[section]
\newtheorem{corollary}[theorem]{Corollary}
\newtheorem{lemma}[theorem]{Lemma}
\newtheorem{proposition}[theorem]{Proposition}

\newtheorem{remark}[theorem]{Remark}
\def\r{\;{\cal R}\;}
\def\l{\;{\cal L}\;}
\def\rs{\;{\cal R^\ast}\;}
\def\ls{\;{\cal L^\ast}\;}

\def\a{{\overline{a}}}
\def\b{{\overline{b}}}
\def\x{{\overline{x}}}
\def\y{{\overline{y}}}
\def\z{{\overline{z}}}

\def\exy{{\overline{xy}}}

\begin{document}
\title{Quasi-ideal transversals of abundant semigroups and spined products}
\author{{\Large Jehan Al-Bar and James Renshaw\footnote{Communicating author}}\\School of Mathematics\\University of Southampton\\Southampton, SO17 1BJ\\England\\Email: j.h.renshaw@maths.soton.ac.uk\\jaal\_bar@hotmail.com}
\date{}
\maketitle

\begin{abstract}
We provide a new and much simpler structure for quasi-ideal adequate transversals of abundant semigroups in terms of spined products, which is similar in nature to that given by Saito for weakly multiplicative inverse transversals of regular semigroups~\cite{saito}. As a consequence we deduce a similar result for multiplicative transversals of abundant semigroups and also consider the case when the semigroups are in fact regular and provide some new structure theorems for inverse transversals.
\end{abstract}

\medskip


\medskip


\section{Introduction}
Much has been written in recent years on adequate transversals of abundant semigroups and the focus has in the main been on determining which of the well-known results and techniques from the theory of inverse transversals carry over to adequate transversals. In particular some effort has been extended to determine structure theorems based on the coordinatisation, $x = (xx^0)x^{00}(x^0x)$ of elements in inverse transversals. We adopt a slightly different approach in this paper and use a structure based on two important subsets of the abundant semigroup, denoted normally by $L$ and $R$. For adequate transversals in general, these do not appear to be subsemigroups, unlike the case for inverse transversals, but it has recently been shown in~\cite{albar-renshaw} that in certain cases they are subsemigroups and in particular when the transversal is a {\em quasi-ideal}. We have therefore attempted to model our structure on a similar one for inverse transversals given by Saito in~\cite{saito}. We refer the reader to~\cite{blyth} for an excellent survey article on inverse transversals and to \cite{fountain1} and \cite{fountain2} for background material on abundant and adequate semigroups.

\bigskip

Unless otherwise stated the terminology and notation for semigroups will be that of~\cite{howie}. Let $S$ be a semigroup and define a left congruence ${\cal R}^\ast$ on $S$ by
$$
{\cal R}^\ast=\{(a,b)\in S\times S\; |\; xa=ya \text{ if and only if } xb=yb \text{ for all } x,y\in S^1\}.
$$
The right congruence ${\cal L}^\ast$ is defined dually. It is easy to show that if $a$ and $b$ are regular elements then $a\rs b$ if and only if $a\r b$. We say that a semigroup is {\em abundant} if each ${\cal R}^\ast-$class and each ${\cal L}^\ast-$class contains an idempotent. An abundant semigroup in which the idempotents commute is called {\em adequate}. It is then clear that regular semigroups are abundant and that inverse semigroups are adequate.
Although morphic images of abundant semigroups need not be abundant, it is clear that isomorphic images are abundant.

\begin{lemma}{\rm(\cite[Corollary 1.2]{fountain2})}\label{e-rs-lemma}
Let $e\in E(S)$ and $a \in S$. Then $e\rs a$ if and only if $ea=a$ and for all $x,y\in S^1$, $xa=ya$ implies $xe=ye$.
\end{lemma}

\begin{lemma}{\rm(\cite[Proposition 1.3]{fountain1})}
A semigroup $S$ is adequate if and only if each ${\cal L}^\ast-$class and each ${\cal R}^\ast-$class contain a unique idempotent and the subsemigroup generated by $E(S)$ is regular.
\end{lemma}

If $S$ is an adequate semigroup and $a\in S$ then we shall denote by $a^\ast$ the unique idempotent in $L_a^\ast$ and by $a^+$ the unique idempotent in $R_a^\ast$.

\medskip

If $S$ is an abundant semigroup and $U$ is an abundant subsemigroup of $S$ then we say that $U$ is a {\em $\ast-$subsemigroup} of $S$ if ${\cal L}^\ast_U = {\cal L}^\ast_S\cap(U\times U), {\cal R}^\ast_U = {\cal R}^\ast_S\cap(U\times U)$. It can be shown that $U$ is a $\ast-$subsemigroup of $S$ if and only if for all $a \in U$ there exist $e,f \in E(U)$ such that $e\in L^\ast_a(S), f\in R^\ast_a(S)$ (see~\cite{el-qallali}). We shall make use of this result below.

\medskip

Let $S$ be an abundant semigroup and $S^0$ be an adequate $\ast-$subsemigroup of $S$. We say that $S^0$ is an {\em adequate transversal} of $S$ if for each $x \in S$ there is a unique $\x\in S^0$ and $e,f\in E$ such that
$$
x = e\x f \hbox{ and such that } e\l\x^+ \hbox{ and } f\r\x^\ast.
$$
It is straightforward to show,~\cite{el-qallali}, that such an $e$ and $f$ are uniquely determined by $x$. Hence we normally denote $e$ by $e_x$, $f$ by $f_x$ and the semilattice of idempotents of $S^0$ by $E^0$. The following are either well-known or are easy to check and are worth mentioning separately:

\begin{enumerate}
\item $e_x\rs x$ and $f_x\ls x$,
\item if $x \in S^0$ then $e_x = x^+\in E^0, \x=x, f_x = x^\ast\in E^0$,
\item $e_\x\l e_x$ and hence $e_\x e_x = e_\x$ and $e_x e_\x = e_x$,
\item $f_\x\r f_x$ and hence $f_\x f_x = f_x$ and $f_x f_\x = f_\x$.
\end{enumerate}
\newcounter{local}
\setcounter{local}{\value{enumi}}

We define
$$
I = \{e_x : x \in S\},\qquad\Lambda=\{f_x : x \in S\}
$$
and note that it is easy to show that $I = \{x\in S: x=e_x\}$ and $\Lambda = \{x\in S: x = f_x\}$.
A number of simple observations regarding $x\in I, y \in \Lambda$ are also worth noting:

\begin{enumerate}
\setcounter{enumi}{\value{local}}
\item $e_x = x, \x=f_x=e_\x$,
\item $e_y=\y=f_\y, f_y=y$.
\end{enumerate}
\setcounter{local}{\value{enumi}}

We also have that for all $x \in S$
\begin{enumerate}
\setcounter{enumi}{\value{local}}
\item $\overline{e_x} = e_\x = \x^+ = f_{e_x}$,
\item $\overline{f_x} = f_\x = \x^\ast = e_{f_x}$,
\item $\x = e_\x xf_\x$.
\end{enumerate}

\smallskip

A number of interesting properties of adequate transversals can be described in terms of the set product $\Lambda I$ and we refer the reader to~\cite{blyth} and~\cite{albar-renshaw} for more details. In particular we say that $S^0$ is a {\em quasi-ideal} of $S$ if $S^0SS^0\subseteq S^0$ or equivalently~\cite[Proposition 2.2]{chen} if $\Lambda I\subseteq S^0$. In~\cite{chen}, Chen proved a structure theorem for such transversals based on $I,\Lambda$ and $S^0$.
We follow a slightly different approach here and define the following two important subsets of $S$ which will play a crucial role in our structure theorem below. Let $S$ be an abundant semigroup with an adequate transversal $S^0$. Define
$$
R = \{x \in S : e_x = e_\x\},\ L = \{x \in S : f_x = f_\x\}
$$
These sets, which in general do not appear to be subsemigroups of $S$, can in fact be described in a number of different ways.
\begin{theorem}{\ \rm(\cite[Theorem 3.1]{albar-renshaw})}
$$R = \{x\in S:e_x \in E^0\} = \{x\in S:x = \x f_x\} = \{x\in S:x = e_\x x\}$$
$$=\{x \in S: x\rs \x\} = \{x\in S:\x = e_x\x\}.$$ 
\end{theorem}

\smallskip

We can summarise the relationship between the various subsets defined with the aid of the following useful diagram
$$
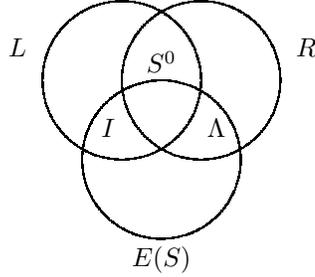
\xy 0;/r10pc/:
{\ellipse<30pt>{}};p;p+(0.25,0),{\ellipse<30pt>{}};
p;p+(-0.125,-.25),{\ellipse<30pt>{}};
p+(0,0.25)*+!D{S^0},+(-0.4,0.05)*+!RD{L},+(0.8,0)*+!LD{R},+(-0.4,-0.55)*+!U{E(S)},+(-.12,0.35)*+!R{I},,+(.24,0)*+!L{\Lambda}
\endxy
$$

\smallskip

It was demonstrated in~\cite[Corollary 3.6]{albar-renshaw} that $LR=S$. However we also have that 
if $x\in R, y \in L$ then $xy = e_\x x y f_\y \in S^0SS^0$ and so if $S^0$ is a quasi-ideal of $S$ then $RL\subseteq S^0$. Conversely if $RL\subseteq S^0$ then in particular $\Lambda I\subseteq S^0$ and so $S^0$ is a quasi-ideal. We have therefore proved
\begin{lemma}
Let $S^0$ be an adequate transversal of an abundant semigroup $S$. Then $S^0$ is a quasi-ideal of $S$ if and only if $RL\subseteq S^0$.
\end{lemma}
\smallskip

The following result will be used extensively.

\begin{lemma}{\ \hbox{\rm(\cite[Theorem 3.4]{albar-renshaw})}}\label{R-S-L-lemma}
Let $x\in R, y \in S, z\in L$. Then
\begin{enumerate}
\item $e_{e_yx}=e_ye_x$,\quad $f_{zf_y}=f_zf_y$,
\item $\overline{e_yx} = e_\y\x$,\quad $\overline{zf_y} = \overline{z}f_\y$,
\item $f_{e_yx}=\left(e_\y\x\right)^\ast f_x$, \quad $e_{zf_y}=e_z\left({\overline z} f_\y\right)^+$.
\end{enumerate}
\end{lemma}

\begin{theorem}{\ \hbox{\rm(Cf. \cite[Theorem 3.12]{albar-renshaw} \& \cite[Proposition 2.1]{luo})}}\label{product-theorem}
Let $S$ be an abundant semigroup with an adequate transversal $S^0$ and let $x,y\in S$. Suppose that $\x f_xe_y, f_xe_y\y\in S^0$. Then
\begin{enumerate}
\item $\overline{xy}=\x f_xe_y\y$,
\item $e_{xy}=e_x(\x f_xe_y)^+$,
\item $f_{xy}=(f_xe_y\y)^\ast f_y$.
\end{enumerate}
\end{theorem}

\section{Quasi-ideals and spined products}

We say that $S$ is {\em left (resp. right) adequate} if $S$ is abundant and every ${\cal R}^\ast-$class (resp. ${\cal L}^\ast-$class) contains a unique idempotent.

\begin{theorem}{\ \rm(\cite[Theorem 3.14]{albar-renshaw})}\label{left-adequate-theorem}
Let $S$ be an abundant semigroup with an adequate transversal $S^0$. The following are equivalent:
\begin{enumerate}
\item $S$ is left adequate;
\item $\Lambda = E^0$;
\item $R = S^0$;
\item $L=S$;
\item $I=E$.
\end{enumerate}
\end{theorem}

The following, and its left-right dual, is then immediate and will be used later without further reference.

\begin{corollary}
Let $S$ be a left adequate semigroup with an adequate transversal $S^0$. Then for all $x \in S, f_x = f_\x$ and so $x = e_x\x$.
\end{corollary}

\begin{theorem}
If $S^0$ is a quasi-ideal adequate transversal of an abundant semigroup $S$ then $L$ and $R$ are subsemigroups of $S$, $L$ is left adequate and $R$ is right adequate with $S^0$ a common quasi-ideal adequate transversal of both $L$ and $R$.
\end{theorem}

\begin{proof}
We deal with $L$ only, the case for $R$ being similar. That $L$ is a subsemigroup of $S$ follows from~\cite[Corollary 4.8]{albar-renshaw}. Let $x\in S^0$ and notice that $x^+{\cal R}^\ast_{S^0} x$ and hence $x^+{\cal R}^\ast_{S} x$ since $S^0$ is a $\ast-$subsemigroup of $S$. It clearly follows then that $x^+{\cal R}^\ast_Lx$. Similarly $x^\ast{\cal L}^\ast_Lx$ and so $S^0$ is a $\ast-$subsemigroup of $L$. For any $x \in L$ we see that $x = e_x\x f_x = e_x\x f_\x$ and $e_x, f_\x\in E(L)$ with $e_x\l\x^+, f_\x\r\x^\ast$. Since $\x$ is clearly unique in this sense then $S^0$ is an adequate transversal of $L$. It is easy to see that $S^0$ is a quasi-ideal transversal of $L$. It also follows from the above observations that $I(L) = I(S) = E(L)$ and so by Theorem~\ref{left-adequate-theorem}, $L$ is left adequate.
\end{proof}

Suppose in what follows that $S^0$ is a quasi-ideal adequate transversal of an abundant semigroup $S$ and $R$ and $L$ are as defined above. Consider the set
$$
T=L|\times|R = \{(x,a) \in L\times R: \x = \a\}
$$
and notice that for $a \in R,y \in L, ay = e_\a ayf_\y \in S^0SS^0\subseteq S^0$.
Also, using Lemma~\ref{R-S-L-lemma}~(2), we see that for $(x,a), (y,b) \in T$
$$
\overline{e_xay} = e_\x\overline{ay} = e_\x ay = e_\a ay = ay,
$$
and in a similar way $\overline{ayf_b} = ay$.  We can therefore define a multiplication on $T$ by
$$
(x,a)(y,b) = (e_xay,ayf_b).
$$
Notice that the set $T$ can be viewed as the spined product, or pullback, within the category of sets, of $L$ and $R$ with respect to $S^0$.
$$
\xymatrix{
L|\times |R\ar[d]\ar[r]&R\ar[d]^{x\mapsto\x}\\
L\ar[r]_{x\mapsto\x}&S^0
}
$$

Although not immediately apparent why, we shall see that the following property of $T$ is very important.
\begin{lemma}\label{product-lemma}
If $(x,a)\in T$ then $xf_a = e_xa$.
\end{lemma}

\begin{proof}
$$
xf_a=e_x\x f_a=e_x\a f_a=e_xa\\
$$
\end{proof}

\begin{lemma}
$T$ is an abundant semigroup and $T\cong S$.
\end{lemma}

\begin{proof}
Consider the well-defined map $\phi : S \to T$ given by
$$
x\phi = (e_x\x,\x f_x).
$$
Notice that from Lemma~\ref{R-S-L-lemma}~(1 and 2), $\overline{e_x\x} = e_\x\x = \x = \x f_\x = \overline{\x f_x}$ and that $e_x = e_xe_\x = e_{e_x\x}$ and $f_x = f_\x f_x = f_{\x f_x}$.
So if $x\phi = y\phi$ then $e_x\x = e_y\y$ and $\x f_x = \y f_y$ from which we deduce that $e_x = e_y, f_x = f_y$ and $\x = \y$. Consequently $x = y$ and we deduce that $\phi$ is one-to-one.

Let $(x,a)\in T$ and let $y = e_xa = xf_a$. Then from Lemma~\ref{R-S-L-lemma}~(1) we see that $e_y = e_xe_a = e_xe_\a = e_xe_\x = e_x$. Also $f_y = f_xf_a = f_\x f_a = f_\a f_a = f_a$. Finally from Lemma~\ref{R-S-L-lemma}~(2), $\y = \overline{e_xa} = e_\x\a = e_\x\x = \x = \a$. Hence $(x,a) = (e_x\x,\a f_a) = (e_y\y,\y f_y) = y\phi$ and so $\phi$ is onto.

From Lemma~\ref{R-S-L-lemma}~(1) and Theorem~\ref{product-theorem}, if $x,y \in S$ then
$$
\begin{array}{rl}
x\phi y\phi&= (e_x\x,\x f_x)(e_y\y,\y f_y)\\
& = (e_{e_x \x}\x f_xe_y\y,\x f_xe_y\y f_{f_y})\\
& = (e_x(\x f_xe_y)^+\x f_xe_y\y,\x f_xe_y\y(f_xe_y\y)^\ast f_y)\\
& = (e_{xy}\exy,\exy f_{xy})\\
&=(xy)\phi
\end{array}
$$
and so $\phi$ is an isomorphism. Hence $T$ is isomorphic to $S$ and is therefore abundant.
\end{proof}

\begin{lemma}\label{regular-lemma}
Let $S$ and $T$ be as above. Then
$$
E(T) = \{(e_y\y,\y f_y) : y \in E(S)\}.
$$
Moreover, if $(x,a)\in E(T)$ then $x,a\in Reg(S)$.
\end{lemma}

\begin{proof}
Suppose $(x,a)\in E(T)$. Then $(x,a)(x,a) = (x,a)$ and so $(e_xax,axf_a) = (x,a)$ which means that $e_xax = x$ and $axf_a = a$. So by Lemma~\ref{product-lemma}, $x = xf_ax$ and $a = ae_xa$ and hence $x,a\in Reg(S)$. Also for some $y \in S, (x,a) = y\phi = (e_y\y,\y f_y)$ and it is easy to see that $y\in E(S)$ since $\phi$ is an isomorphism. Hence $x=e_y\y, a=\y f_y$ with $y\in E(S)$.
\end{proof}

\bigskip

Notice that if $x\in S$ then $\x\phi = (e_\x\x,\x f_\x) = (\x,\x)$. So let
$$
T^0 = \{(\x,\x) \in L\times R : x\in S\}.
$$
Since $\phi$ is an isomorphism and since $S^0$ is a quasi-ideal adequate transversal of $S$ then $T^0$ is a quasi-ideal adequate transversal of $T$.
Consequently we have shown

\begin{proposition}\label{converse-proposition}
Let $S^0$ be a quasi-ideal adequate transversal of an abundant semigroup $S$. Then there exists a left adequate semigroup $L$ and a right adequate semigroup $R$ with $S\cong L|\times|R = T$ and such that $T$ contains a quasi-ideal adequate transversal $T^0\cong S^0$.
\end{proposition}

Notice that the multiplication in $T$, $(x,a)(y,b) = (e_xay,ayf_b)$, relies on the existence of a map $R\times L\to S^0$, given by $(a,x)\mapsto ax$. Moreover we have frequently used some basic properties of $L$ and $R$ such as $e_a(ax) = ax$. In order to consider a converse for Proposition~\ref{converse-proposition} suppose that $L$ is a left adequate semigroup and $R$ is a right adequate semigroup and suppose that $S^0$ is a common quasi-ideal adequate transversal of $L$ and $R$ with semilattice of idempotent $E^0$. Then for each $x\in L, a\in R$ we have factorisations
$$
x = e_x\x f_x,\quad a = e_a\a f_a
$$
where $\x,\a\in S^0$, $e_x\in I(L)=E(L), f_x\in \Lambda(L)=E^0$ and $e_a\in I(R)=E^0, f_a\in \Lambda(R)=E(R)$. Notice that if $y\in S^0$ then $e_y, f_y\in E^0$ and have the same value whether we consider $S^0$ to be a transversal of $L$ or of $R$. The notation $e_x$ and $e_a$, for $x\in L, a \in R$ should not cause any confusion.

\begin{theorem}\label{main-theorem}
Let $S^0$ be an adequate semigroup with semilattice of idempotents $E^0$. Let $L$ be a left adequate semigroup and $R$ a right adequate semigroup and suppose that $S^0$ is a common quasi-ideal adequate transversal of $L$ and $R$. Let $R\times L\to S^0$ be a map denoted by $(a,x)\mapsto a\ast x$ and suppose that $\ast$ satisfies
\begin{enumerate}
\item for all $y,z\in L, a,b\in R$ with $\y=\b$,
$(a\ast y)f_b\ast z = a\ast e_y(b\ast z)$;
\item if $a\in S^0$ or $x\in S^0$ then $a\ast x = ax$.
\end{enumerate}
Let $T = \{(x,a)\in L\times R : \x = \a\}$ and define
$$
(x,a)(y,b)=(e_x(a\ast y),(a\ast y)f_b).
$$
Then $T$ is an abundant semigroup with a quasi-ideal adequate transversal $T^0$ with $T^0\cong S^0$. Moreover every quasi-ideal adequate transversal can be constructed in this way.
\end{theorem}

\begin{proof}
Notice that property (1) is essentially the analogue of the important property given in Lemma~\ref{product-lemma}, property (2) says that $\ast$ is an extension of the multiplication in $L$ and $R$.
We first prove a few lemmas that will be useful later. In what follows, $R, L,S^0$ and $T$ are as stated in the Theorem.

The following result is needed to show that $T$ is abundant.
\begin{lemma}\label{property3}
For all $a,b \in R, x,y \in L$, if $a\ast x = b\ast x$ then $a\ast e_x = b\ast e_x$ and if $a\ast x = a\ast y$ then $f_a\ast x = f_a\ast y$.
\end{lemma}
\begin{proof}
Suppose that $a\ast x = b\ast x$. Then we have the following argument
$$
\begin{array}{rclr}
a\ast x&=&b\ast x&\\
a\ast e_x\x&=&b\ast e_x\x&x\in L\\
a\ast e_x(e_\x\ast\x)&=&b\ast e_x(e_\x\ast\x)&\\
a\ast e_{e_x}(e_\x\ast\x)&=&b\ast e_{e_x}(e_\x\ast\x)&e_x = e_{e_x}\\
(a\ast e_x)f_{e_\x}\ast \x&=&(b\ast e_x)f_{e_\x}\ast \x&\hbox{\rm property (1)}\\
(a\ast e_x)f_{e_\x}\x&=&(b\ast e_x)f_{e_\x}\x&\x\in S^0\\
(a\ast e_x)f_{e_\x}\x^+&=&(b\ast e_x)f_{e_\x}\x^+&\x\rs\x^+\\
(a\ast e_x)f_{e_\x}\ast \x^+&=&(b\ast e_x)f_{e_\x}\ast \x^+&\\
a\ast e_{e_x}(e_\x\ast\x^+)&=&b\ast e_{e_x}(e_\x\ast\x^+)&\\
a\ast e_x\x^+&=&b\ast e_x\x^+&e_\x=\x^+\\
a\ast e_x&=&b\ast e_x&e_x\l\x^+\\
\end{array}
$$
The other identity follows in a similar manner.
\end{proof}

\begin{lemma}\label{sub-lemma}
If $a\in R, x \in L, g\in E(R), h\in E(L)$ with $\overline{g}=f_\a, \overline{h}=e_\x$ then $ag\ast hx = e_a(ag\ast hx)f_x$. In particular $a\ast x = e_a(a\ast x)f_x$.
\end{lemma}
\begin{proof}
Notice first that $e_g = \overline{g}, f_g = g$ and $e_h = h,f_h = \overline{h}$. Since $S^0$ is a quasi-ideal transversal of $R$ then by Theorem~\ref{product-theorem} and since $e\in \Lambda, f\in I$
$$
\overline{ag}=\a f_ae_g\overline{g}=\a f_a \overline{g}=\a f_af_\a=\a f_\a = \a,
$$
and
$$
\overline{hx}=\overline{h} f_he_x\x=\overline{h}e_x\x=e_\x e_x\x=e_\x\x = \x.
$$
Also
$$
f_{ag}=\left(f_ae_g\overline{g}\right)^\ast f_g=\left(f_a\overline{g}\right)^\ast g=\left(f_a f_\a\right)^\ast g=f_\a g=\overline{g}g=g,
$$
and
$$
e_{hx}=e_h\left(\overline{h}f_he_x\right)^+=h\left(e_\x e_x\right)^+=h\left(e_\x\right)^+=he_\x=h\overline{h}=h.
$$
So using properties (1) and (2) we have
$$
\begin{array}{rlr}
e_a(ag\ast hx)&=e_\a(ag\ast hx)&\\
&=e_\a\ast e_a(ag\ast hx)&\\
&=(e_\a\ast a)f_{ag}\ast hx&\hbox{\rm since }\overline{a}=\overline{ag}\\
&=ag\ast hx&\hbox{\rm since }f_{ag}=g.
\end{array}
$$
In a similar manner
$$
\begin{array}{rl}
(ag\ast hx)f_x&=(ag\ast hx)f_\x\\
&=(ag\ast hx)f_x\ast f_\x\\
&=ag\ast(e_{hx}(x\ast f_\x))\\
&=ag\ast hx
\end{array}
$$
\end{proof}

\begin{lemma}\label{e-f-lemma}
Let $a,b \in R, y,z\in L$ and suppose that $\y = \b$. Then $e_{a\ast y} = e_{(a\ast y)f_b}$ and $f_{e_y(b\ast z)} = f_{b\ast z}$.
\end{lemma}
\begin{proof}
Using Lemma~\ref{R-S-L-lemma}~(3) and Lemma~\ref{sub-lemma} we see that
$$
\begin{array}{ll}
e_{(a\ast y)f_b}&= e_{a\ast y}\left(\overline{a\ast y}f_\b\right)^+ = e_{a\ast y}\left(\overline{a\ast y}f_\y\right)^+\\
&=e_{a\ast y}\left((a\ast y)f_y\right)^+ = e_{a\ast y}\left(a\ast y\right)^+\\
&=e_{a\ast y}e_{a\ast y} = e_{a\ast y}\\
\end{array}
$$
and 
$$
\begin{array}{ll}
f_{e_y(b\ast z)}&= \left(e_\y\overline{(b\ast z)}\right)^\ast f_{b\ast z}= \left(e_\b\overline{(b\ast z)}\right)^\ast f_{b\ast z}\\
&=\left(e_b(b\ast z)\right)^\ast f_{b\ast z}=\left(b\ast z\right)^\ast f_{b\ast z}\\
&=f_{b\ast z}f_{b\ast z}\ = f_{b\ast z}\\
\end{array}
$$

\end{proof}

Now to show that $T$ is a semigroup we use property (1) together with Lemma~\ref{R-S-L-lemma}~(1), Lemma~\ref{sub-lemma} and Lemma~\ref{e-f-lemma} as follows:
$$
\begin{array}{rl}
\left((x,a)(y,b)\right)(z,c)&=(e_x(a\ast y),(a\ast y)f_b)(z,c)\\
&=(e_{e_x(a\ast y)}((a\ast y)f_b\ast z),((a\ast y)f_b\ast z)f_c)\\
&=(e_xe_{a\ast y}((a\ast y)f_b\ast z),((a\ast y)f_b\ast z)f_c)\\
&=(e_xe_{(a\ast y)f_b}((a\ast y)f_b\ast z),((a\ast y)f_b\ast z)f_c)\\
&=(e_x(a\ast (e_y(b\ast z))),(a\ast(e_y(b\ast z)))f_{e_y(b\ast z)}f_c)\\
&=(e_x(a\ast (e_y(b\ast z))),(a\ast(e_y(b\ast z)))f_{(b\ast z)}f_c)\\
&=(e_x(a\ast (e_y(b\ast z))),(a\ast(e_y(b\ast z)))f_{(b\ast z)f_c})\\
&=(x,a)(e_y(b\ast z),(b\ast z)f_c)\\
&=(x,a)\left((y,b)(z,c)\right).\\
\end{array}
$$
To show that $T$ is abundant we require the following lemmas.
\begin{lemma}\label{sub-R-S-L-lemma}
For $a\in R, x \in L, (a\ast e_x){\cal R}^\ast_L(a\ast x){\cal L}^\ast_L(f_a\ast x)$.
\end{lemma}
\begin{proof}
First notice that using Lemma~\ref{sub-lemma} we have
$$
\begin{array}{rlcrl}
(a\ast e_x)x &= (a\ast e_x)\ast x&\ {\rm and}\ &a(f_a\ast x) &= a\ast(f_a\ast x)\\
&=(a\ast e_x)f_{e_\x}\ast x&&&=a\ast e_{f_\a}(f_a\ast x)\\
&=a\ast e_x(e_\x\ast x)&&&=(a\ast f_\a)f_a\ast x\\
&=a\ast x&&&=a\ast x.
\end{array}
$$
So we can immediately deduce that if $z,w \in L$ are such that $z(a\ast e_x) = w(a\ast e_x)$ then $z(a\ast x) = w(a\ast x)$ and if $c,d\in L$ are such that $(f_a\ast x)c = (f_a\ast x)d$ then $(a\ast x)c = (a\ast x)d$.
In addition
$$
\begin{array}{rlcrl}
z(a\ast x)&=z\ast(a\ast x)&\ {\rm and }\ &(a\ast x)c&=(a\ast x)\ast c\\
&=z\ast e_{\a}(a\ast x)&&&=(a\ast x)f_\x\ast c\\
&=(z\ast\a)f_a\ast x&&&=a\ast e_x(\x\ast c).
\end{array}
$$
Hence, if conversely we have $z(a\ast x)=w(a\ast x)$ then $(z\ast\a)f_a\ast x = (w\ast\a)f_a\ast x$ and so from Lemma~\ref{property3} we see that $(z\ast\a)f_a\ast e_x = (w\ast\a)f_a\ast e_x$ and therefore $z(a\ast e_x)=w(a\ast e_x)$. In a similar way if
$(a\ast x)c = (a\ast x)d$ then $(f_a\ast x)c = (f_a\ast x)d$. The result then follows.
\end{proof}

\begin{lemma}\label{ET-lemma}
$E(T) = \{(x,a) \in T: \x = a\ast x\}$.
\end{lemma}
\begin{proof}
First notice that $(x,a)(x,a) = (e_x(a\ast x),(a\ast x)f_a)$. So if $a\ast x = \x$ then $(x,a)^2 = (e_x\x,\x f_a) = (e_x\x f_x,e_a\a f_a) = (x,a)$. Conversely, if $(x,a)^2 = (x,a)$ then $e_x(a\ast x) = x$ and so
$$
\begin{array}{ll}
\x&=e_\x x = e_\x\ast x = e_\a\ast x = e_\a\ast e_x(a\ast x)\\
&=(e_\a\ast x)f_a\ast x = (e_\x e_x\x)f_a\ast x = \a f_a\ast x\\
&= a\ast x
\end{array}
$$
\end{proof}

Now let $(x,a) \in T$ and notice that by Lemma~\ref{ET-lemma}, $(e_x,e_\x), (f_\a,f_a) \in E(T)$.
Now
$$
\begin{array}{rl}
(e_x,e_\x)(x,a)&=(e_{e_x}(e_\x\ast x),(e_\x\ast x)f_a)\\
&=(e_x(e_\x\ast x),(e_\x\ast x)f_a)\\
&=(e_x\x f_x,\x f_xf_a)\\
&=(x,\x f_\x f_a)\\
&=(x,\a f_a)\\
&=(x,a).
\end{array}
$$
Suppose then that
$$
(y,b)(x,a) = (z,c)(x,a).
$$
Then
$$
(e_y(b\ast x),(b\ast x)f_a) = (e_z(c\ast x),(c\ast x)f_a)
$$
and so $(b\ast x)f_a f_\a f_x= (c\ast x) f_af_\a f_x$ which means, by Lemma~\ref{sub-lemma} that $b\ast x = c\ast x$. Also $e_y(b\ast x) = e_z(b\ast x)$ and so from Lemma~\ref{sub-R-S-L-lemma}, $e_y(b\ast e_x)=e_z(b\ast e_x)$. In addition it follows from Lemma~\ref{property3} that $b\ast e_x=c\ast e_x$ and so $e_y(b\ast e_x)=e_z(c\ast e_x)$ while $(b\ast e_x)f_{e_\x} = (c\ast e_x)f_{e_\x}$ as required to prove
$$
(y,b)(e_x,e_\x) = (z,c)(e_x,e_\x)
$$
and so $(x,a){\cal R}^\ast_T (e_x,e_\x)$. Similarly $(x,a){\cal L}^\ast_T(f_\a,f_a)$ and so $T$ is abundant.

\medskip

Let $T^0 = \{(s,s) : s\in S^0\}$. Then $T^0$ is a subsemigroup of $T$ since if $s,t \in S^0$ then
$$
\begin{array}{rl}
(s,s)(t,t)&=(e_s(s\ast t),(s\ast t)f_t)\\
&=(s\;t,s\;t)
\end{array}
$$
and clearly $T^0\cong S^0$. Notice also that $E(T^0) = \{(s,s) : s\in E^0\}$ and since $E^0$ is a semilattice then so is $E(T^0)$. Hence $T^0$ is an adequate subsemigroup of $T$. It is also a $\ast-$subsemigroup of $T$ as, from above we see that for all $s\in S^0$, $(e_s,e_s){\cal R}^\ast_T(s,s)$ and $(f_s,f_s){\cal L}^\ast_T(s,s)$ and that $(e_s,e_s) = (s,s)^+, (f_s,f_s) = (s,s)^\ast$.

Suppose then that $x\in L, a \in R$ and consider
$$
\begin{array}{rl}
(e_x,e_\x)(\x,\x)(f_\a,f_a)&=(e_x(e_\x\ast\x),(e_\x\ast\x)f_\x)(f_\a,f_a)\\
&=(e_x\x,\x)(f_\a,f_a)\\
&=(x,\x)(f_\x,f_a)\\
&=(e_x(\x\ast f_\x),(\x\ast f_\x)f_a)\\
&=(e_x\x,\a f_a)\\
&=(x,a).
\end{array}
$$
To demonstrate that $(\x,\x)$ is unique with respect to these properties, suppose that
$$
(z,b)(\y,\y)(w,c)=(x,a)
$$
with $(z,b)\l(\y,\y)^+$ and $(w,c)\r(\y,\y)^\ast$ with $(z,b),(w,c)\in E(T)$. Then
$$
(z,b)=(z,b)(e_\y,e_\y)
=(e_z(b\ast e_\y),(b\ast e_\y)f_{e_\y})
=(e_zbe_\y,be_\y))
$$
and
$$
(e_\y,e_\y)=(e_\y,e_\y)(z,b)
=(e_\y(e_\y\ast z),(e_\y\ast z)f_b)
=(e_\y z,e_\y zf_b)
$$
So $z = e_zbe_\y$ and therefore $ze_\y=z$. Also, $e_\y z = e_\y$ and so $z\l\y^+$. Notice that $b = b\ast e_\y\in S^0$. Since $(z,b)\in E(T)$ then $z = e_z(b\ast z)=e_zbz=z^2$, since $z=e_zb$.
In a similar way we can deduce that
$$
w = f_\y w, c = wf_c, f_\y = cf_\y, e_wf_\y=f_\y
$$
and so see that $w\r\y^\ast$ and $w\in E(L)$. Finally
$$
\begin{array}{rl}
(z,b)(\y,\y)(w,c)&=(e_zb\y,b\y)(w,c)=(z\y,b\y)(w,c)\\
&=(e_{z\y}b\y w,b\y wf_c)=(e_{e_zb\y}b\y w,b\y wf_c)\\
&=(e_ze_{b\y}b\y w,b\y c)=(e_zb\y w,b\y c)\\
&=(z\y w,b\y c)
\end{array}
$$
and so $x = z\y w$ and hence $\y = \x$ by the uniqueness of $\x$. Consequently, $T^0$ is an adequate transversal of $T$. Finally
$$
\begin{array}{rl}
(\y,\y)(x,a)(\z,\z)&=(e_\y\y x,\y xf_a)(\z,\z)\\
&=(\y x,\y xf_a)(\z,\z)\\
&=(e_{\y x}\y xf_a\z,\y xf_a\z f_\z)\\
&=(\y xf_a\z,\y xf_a\z)
\end{array}
$$
and since $\y xf_a\z \in S^0SS^0 = S^0$ then $T^0$ is a quasi-ideal of $T$.

\bigskip

Conversely let $S^0$ be a quasi-ideal adequate transversal of an abundant semigroup $S$ and for all $a\in R, x \in L$ let $a\ast x = ax$. Most of the converse follows from Proposition~\ref{converse-proposition} and so we need only check that $\ast$ satisfies the given conditions.
\begin{enumerate}
\item Suppose that $y,z\in L, a,b\in R$ and $\y=\b$. Then
$$
((a\ast y)f_b\ast z)=ayf_bz
=ae_y\y f_bz
=ae_y\b f_bz
=ae_ybz
=a\ast (e_y(b\ast z)).
$$
\item It is clear that if $a\in S^0$ or $x\in S^0$ then $a\ast x = ax$;
\end{enumerate}
The proof of Theorem~\ref{main-theorem} is now complete.
\end{proof}

\begin{remark}\label{main-remark}
\end{remark}
Let $T$ be as in Theorem~\ref{main-theorem}. Notice that $(x,a)\in R(T)$ if and only if $e_{(x,a)} = e_{\overline{(x,a)}}$ or in other words $(e_x,e_\x) = (e_\x,e_\x)$ which is clearly equivalent to $x\in R(L) = S^0$ and so $x=\a$. Define a map $R\to T$ by $a\mapsto (\a,a)$ and notice that $(\a,a)(\b,b) = (e_\a(a\ast\b),(a\ast\b)f_b)=(a\b,ab)$. From Theorem~\ref{product-theorem} we can easily deduce that $a\b = \overline{ab}$ and hence $R(T)\cong R$. 
Also, we see from the above proof that $\Lambda(T) = \{(f_\a,f_a) : a\in R\} = \{(\a,a) : a\in \Lambda(R)\}$.
In a similar way $L(T)\cong L$ and $I(T) = \{(x,\x) : x \in I(L)\}$.

\bigskip
\def\tensor{\otimes}

We have consequently reduced the problem of finding the structure of abundant semigroups with quasi-ideal adequate transversals to one of determining the structure of quasi-ideal adequate transversals of left adequate and right adequate semigroups. For this we can make use of the alternative description of abundant semigroups with quasi-ideal adequate transversals given by Chen in \cite{chen}. Suppose then that $S^0$ is an adequate semigroup with a semilattice of idempotents $E^0$. Let $A$ be a set disjoint from $S^0$ and let $I=E^0\cup A$. Suppose also that $I$ is a right $E^0-$act. In other words there is a map $I\times E^0\to I$ given by $(x,s)\mapsto x\tensor s$ such that $(x\tensor s)\tensor t = x\tensor (st)$. Suppose also that there is a right $E^0-$map $\phi : I\to E^0$ such that $\phi|_{E^0} = 1_{E^0}$ and such that $x\tensor x\phi = x$ and $(x\tensor s)\phi = x\phi\tensor s = x\phi s$.

\begin{theorem}{\ \rm(Cf. \cite[Theorem 3.1]{chen})}\label{chen-theorem}
Let $S^0$ be an adequate semigroup with semilattice of idempotent $E^0$ and suppose that $I, \tensor$ and $\phi$ are as above. Suppose there is a map $E^0\times I \to S^0$ given by $(f,e)\mapsto f\ast e$ which satisfies
\begin{enumerate}
\item for all $a,b\in E^0$, $a(f\ast e)b = af\ast(e\tensor b)$,
\item for all $f,g\in E^0$, $f\ast g = fg$
\item for all $e\in I$, $e\phi\ast e=e\phi$.
\end{enumerate}
Let
$$
T = \{(e,x)\in I\times S^0: e\phi=x^+\}
$$
and define a multiplication on $T$ by
$$
(e,x)(g,w) = (e\tensor a^+,a)
$$
where $a = x(x^\ast\ast g)w\in S^0$. Then $T$ is a left adequate semigroup with a quasi-ideal adequate transversal $T^0=\{(a^+,a) : a\in S^0\}$ with $T^0$ isomorphic to $S^0$. Moreover every quasi-ideal transversal of a left adequate semigroup can be constructed in this way.
\end{theorem}

Notice that this is slightly different from that quoted in~\cite[Theorem 3.1]{chen}. First, as $\Lambda = E^0$, we have replaced the notation $[f,e]$ by $f\ast e$ in keeping with the notation we have used in this paper; and secondly not all of Chen's conditions (B5) are actually needed for the proof as the part relating to the regularity of $[f\psi,e]$ is only used in the proof of the converse of his Lemma 3.7. However we can alternatively use his property (B4) to deduce that $a_1^+ = a_2^+$ and so deduce that $e_1=e_2$ (infact $a_1^+ = [a_2^+,e_2\tensor a_1^+] = [a_2^+,e_2]a_1^+ = [a_2^+,e_1\tensor a_2^+]a_1^+ = [a_2^+,e_1]a_2^+a_1^+ = a_1^+a_2^+$ and similarly $a_2^+ = a_2^+a_1^+$). We leave the rest of the details to the interested reader.

\section{Multiplicative transversals}

An adequate transveral $S^0$ of an abundant semigroup $S$ is said to be {\em multiplicative} if $\Lambda I \subseteq E^0$ and is said to be {\em weakly multiplicative} if $\overline{\Lambda I}\subseteq E^0$ where $\overline{\Lambda I} = \{\overline{li} : l\in \Lambda, i \in I\}$. It was shown in~\cite[Theorem 5.3]{albar-renshaw} that $S^0$ is multiplicative if and only if $S^0$ is weakly multiplicative and $S^0$ is a quasi-ideal of $S$. Consequently we can use Theorem~\ref{main-theorem} to deduce

\begin{theorem}\label{multiplicative-theorem}
Let $S^0$ be an adequate semigroup with semilattice of idempotents $E^0$. Let $L$ be a left adequate semigroup and $R$ a right adequate semigroup and suppose that $S^0$ is a common quasi-ideal adequate transversal of $L$ and $R$. Let $R\times L\to S^0$ be a map denoted by $(a,x)\mapsto a\ast x$ and suppose that $\ast$ satisfies
\begin{enumerate}
\item for all $y,z\in L, a,b\in R$ with $\y=\b$,
$(a\ast y)f_b\ast z = a\ast e_y(b\ast z)$;
\item if $a\in S^0$ or $x\in S^0$ then $a\ast x = ax$;
\item for all $e\in E(R), f \in E(L), e\ast f \in E^0$.
\end{enumerate}
Let $T = \{(x,a)\in L\times R : \x = \a\}$ and define
$$
(x,a)(y,b)=(e_x(a\ast y),(a\ast y)f_b).
$$
Then $T$ is an abundant semigroup with a multiplicative adequate transversal $T^0$ with $T^0\cong S^0$. Moreover every multiplicative adequate transversal can be constructed in this way.
\end{theorem}

\begin{proof}
If $S^0, L, R$ and $T$ are as above then from Theorem~\ref{main-theorem}, $T$ is an abundant semigroup with a quasi-ideal adequate transversal $T^0$ isomorphic to $S^0$. If in addition $(x,a)\in \Lambda(T), (y,b)\in I(T)$ then from Remark~\ref{main-remark} we see that $a\in E(R), x = \a, y \in E(L), b = \y$.
$$
(x,a)(y,b) = (e_\x(a\ast y),(a\ast y)f_\b) = (a\ast y, a\ast y) \in E(T^0)
$$
by property (3).
Consequently $T^0$ is
multiplicative.

\smallskip

Conversely if $S^0$ is a multiplicative transversal of an abundant semigroup $S$ then from the proof of Theorem~\ref{main-theorem} we see that $T=\{(x,a)\in L\times R : \x = \a\}\cong S$ is an abundant semigroup with a quasi-ideal transversal $T^0 = \{(\x,\x) \in L\times R : x\in S\}\cong S^0$ such that properties (1) and (2) hold when $a\ast x = ax$ for $a\in R, x\in L$. Suppose then that $e\in E(R), f \in E(L)$ so that $e\ast f = ef \in \Lambda I \subseteq E^0$ and so property (3) holds as well.
\end{proof}

Notice that although we only specify that $S^0$ is a {\em quasi-ideal transversal} of the semigroups $L$ and $R$, property (3) guarantees that $S^0$ will in fact be a {\em multiplicative transversal} of both $L$ and $R$. To see this, let $\lambda\in \Lambda(R)
, i \in I(R)
$ and notice that
$(\overline{\lambda},\lambda) \in \Lambda(T)$ and $(i,\overline{i}) \in I(T)$ by Remark~\ref{main-remark}. Moreover since $T^0$ is a multiplicative transversal of $T$ then $(\lambda i,\lambda i) = (\overline{\lambda},\lambda)(i,\overline{i}) \in E^0(T)$ and so $\lambda i \in E^0$ and hence $S^0$ is a multiplicative transversal of $R$. In a similar way $S^0$ is a multiplicative transversal of $L$.

\section{The regular case}
Suppose that $S$ is an abundant semigroup with an adequate transversal $S^0$ and suppose that $x\in Reg(S)$, the set of regular elements of $S$. Using the fact that $x\r e_x$ and $x\l f_x$ then from \cite[Theorem 2.3.4]{howie} there exists a unique $x^0\in V(x)$ with $xx^0=e_x$ and $x^0x=f_x$.
In what follows we shall write $x^{00}$ for $\left(x^0\right)^0$.

\begin{theorem}\ {\rm\cite[Theorems 2.3 \& 2.4]{albar-renshaw}}\label{first-regular-theorem}
If $x\in Reg(S)$ then $|V(x)\cap S^0|=1$. Moreover $x^0\in S^0, \x = x^{00}$ and $x^0 = x^{000}$.
Also for all $x \in S, e_x = xx^0, f_x = x^0x$ and
$$I = \{x\in Reg(S) : x = xx^0\} = \{xx^0 : x\in Reg(S)\}$$
and
$$\Lambda = \{x\in Reg(S) : x = x^0x\} = \{x^0x : x \in Reg(S)\}.$$
\end{theorem}

Consequently if $S^0$ is an inverse transversal of a regular semigroup $S$ then $S^0$ is an adequate transversal of the abundant semigroup $S$ and for all $x \in S, e_x = xx^0, f_x = x^0x$ and $\x = x^{00}$. Notice also that in this case for $x,y \in S, \x = \y$ if and only if $x^0 = y^0$ and if $i \in I$ then $\overline{i} = f_i = i^0$.

\medskip

Let us suppose now that all our semigroups are regular and so in particular we assume that $L$ is a left inverse semigroup and $R$ is a right inverse semigroup with a common inverse transversal $S^0$. Suppose also that there is a map $R\times L\to S^0$ given by $(a,x)\mapsto a\ast x$.

\medskip

From Theorem~\ref{main-theorem} we deduce
\begin{theorem}\label{regular-theorem}
Let $S^0$ be an inverse semigroup with semilattice of idempotents $E^0$. Let $L$ be a left inverse semigroup and $R$ a right inverse semigroup and suppose that $S^0$ is a quasi-ideal inverse transversal of both $L$ and $R$. Let $R\times L\to S^0$ be a map denoted by $(a,x)\mapsto a\ast x$ and suppose that $\ast$ satisfies
\begin{enumerate}
\item for all $y,z\in L, a,b\in R$ with $y^0=b^0$,
$(a\ast y)(b^0b)\ast z = a\ast (yy^0)(b\ast z)$;
\item if $a\in S^0$ or $x\in S^0$ then $a\ast x = ax$.
\end{enumerate}
Let $T = \{(x,a)\in L\times R : x^0 = a^0\}$ and define
$$
(x,a)(y,b)=(xx^0(a\ast y),(a\ast y)b^0b).
$$
Then $T$ is a regular semigroup with a quasi-ideal inverse transversal $T^0$ with $T^0\cong S^0$. Moreover every quasi-ideal inverse transversal of a regular semigroup can be constructed in this way.
\end{theorem}

\begin{proof}
If $S^0$ is an inverse semigroup and $L$ and $R$ are as stated in the theorem, then from Theorem~\ref{main-theorem} we see that $T= \{(x,a)\in L\times R: \x = \a\}$ is an abundant semigroup with a quasi-ideal adequate transversal isomorphic to $S^0$. However, from Theorem~\ref{first-regular-theorem} we easily deduce that $x^0 = a^0$ if and only if $\x = \a$.
It is also a fairly simple matter to check that for all $(x,a)\in T, (x,a)(a^0,x^0)(x,a) = (x,a)$ and so $T$ is regular and hence $T^0$ is inverse.

\medskip

The converse follows easily from Theorem~\ref{main-theorem} and Theorem~\ref{first-regular-theorem}.
\end{proof}

In a similar way we can also deduce from Theorem~\ref{chen-theorem}

\begin{theorem}
Let $S^0$ be an inverse semigroup with semilattice of idempotent $E^0$ and suppose that $I$ is a left normal band with semilattice transversal $E^0$.

Let
$$
T = \{(e,x)\in I\times S^0: e^0=xx^{-1}\}
$$
and define a multiplication on $T$ by
$$
(e,x)(g,w) = (e(aa^{-1}),a)
$$
where $a = x((x^{-1}x)g)w\in S^0$. Then $T$ is a left inverse semigroup with a quasi-ideal inverse transversal $T^0=\{(aa^{-1},a) : a\in S^0\}$ with $T^0$ isomorphic to $S^0$. Moreover every quasi-ideal inverse transversal of a left inverse semigroup can be constructed in this way.
\end{theorem}

\begin{proof}
First notice that if $s \in E^0, g \in I$ then since $I$ is left normal we have $sg = sgg^0g = sg^0gg = sg^0\in E^0$.
If $S^0$ and $I$ are as given in the theorem then first define $\phi : I \to E^0$ by $x\phi = x^0$ and the map $I\times E^0\to I$ by $(x,s)\mapsto x\tensor s = xs$ and $E^0\times I\to E^0$ by $s\ast x = sx$. Then it is straightforward to show that $\phi,\ast$ and $\tensor$ satisfy the required properties. Hence by Theorem~\ref{chen-theorem} we can construct $T$ and $T^0$ as in the statement of the theorem and we see that $T^0$ is a quasi-ideal adequate transversal of the left adequate semigroup $T$. That $T$ is regular follows from the observation that for all $(e,x)\in T, (x^{-1}x,x^{-1})\in T$ and $(e,x)(x^{-1}x,x^{-1})(e,x) = (e,x)$.

\medskip

Conversely, if $T^0$ is a quasi-ideal inverse transversal of a left inverse semigroup $T$, then $T^0$ is adequate and $T$ is left adequate and so from Theorem~\ref{chen-theorem} we see that there exist $S^0, E^0, I, \tensor,\phi$ and $\ast$ as in Theorem~\ref{chen-theorem} such that $T\cong \{(e,x) \in I\times S^0 : e\phi = x^+\}$ and $T^0\cong S^0$ with multiplication as given. From the proof of Theorem~\ref{chen-theorem} together with the remarks after Theorem~\ref{first-regular-theorem} we see that $I = I(T)$ and for all $x\in I, s \in E^0, x\tensor s = xs, x\phi = \x = x^0$ and $s\ast x = sx$. In addition we know that $I(T)$ is a left normal band (\cite[Proposition 1.7]{mcalister}). Finally since for all $x \in S^0, x^+ = xx^{-1}, x^\ast = x^{-1}x$ then the result follows.
\end{proof}

Finally from Theorem~\ref{multiplicative-theorem} we can deduce
\begin{theorem}
Let $S^0$ be an inverse semigroup with semilattice of idempotents $E^0$. Let $L$ be a left inverse semigroup and $R$ a right inverse semigroup and suppose that $S^0$ is a common quasi-ideal inverse transversal of $L$ and $R$. Let $R\times L\to S^0$ be a map denoted by $(a,x)\mapsto a\ast x$ and suppose that $\ast$ satisfies
\begin{enumerate}
\item for all $y,z\in L, a,b\in R$ with $y^0=b^0$,
$(a\ast y)(b^0b)\ast z = a\ast (yy^0)(b\ast z)$;
\item if $a\in S^0$ or $x\in S^0$ then $a\ast x = ax$;
\item for all $e\in E(R), f \in E(L), e\ast f \in E^0$.
\end{enumerate}
Let $T = \{(x,a)\in L\times R : x^0 = a^0\}$ and define
$$
(x,a)(y,b)=((xx^0)(a\ast y),(a\ast y)(b^0b)).
$$
Then $T$ is a regular semigroup with a multiplicative inverse transversal $T^0$ with $T^0\cong S^0$. Moreover every multiplicative inverse transversal can be constructed in this way.
\end{theorem}

\begin{proof}
Let $S^0, E^0, L$ and $R$ be as stated in the Theorem. Then by Theorem~\ref{multiplicative-theorem} $T$ is an abundant semigroup with a multiplicative adequate transversal $T^0\cong S^0$. That $T$ is regular follows from the observation that $(x,a)(a^0,x^0)(x,a) = (x,a)$.

The converse is straightforward.
\end{proof}

The authors would like to thank the anonymous referee for pointing out a more efficient proof of Proposition~\ref{converse-proposition}.

\end{document}